\documentclass[12pt]{article}
\usepackage{amssymb , amsmath, amsthm}
\def\picta{\unitlength=1cm
 \begin{picture}(0,0)
\put(4,0){\vector(1,1){2}}
 \put(4, 2){\line(1,-1){.8}}
 \put(5.1, .85){\vector(1,-1){.85}}
 \put(4.7, -.5){ -1}
 \put(8,0){\line(1,1){.8}}
\put(9.2,1.2){\vector(1,1){.8}}
\put(8,2){\vector(1,-1){2}}
\put(8.7, -.5){ +1}o
 \end{picture}\hspace{2cm}}
\begin{document}
\newtheorem{thm}{Theorem}
\newtheorem{cor}{Corollary}
\newtheorem{lem}{Lemma}
\newtheorem{slem}{Sublemma}
\newtheorem{prop}{Proposition}
\newtheorem{defn}{Definition}
%\makeatother
%\def\email#1{\def\@email{#1}}
\title{Some properties of simple minimal knots}
\author{Marc Soret and Marina Ville \\
{\it  \footnotesize{marc.soret@ lmpt.univ-tours.fr \& Marina.Ville@ lmpt.univ-tours.fr}}}
 \date{ }
\maketitle
 \begin{abstract}
 A minimal knot is the intersection of a topologically embedded branched minimal disk in $\mathbb{R}^4$
  $\mathbb{C}^2 $  with a small sphere centered at the branch point. When the lowest order terms in each  coordinate  component    of the embedding of the disk  in $\mathbb{C}^2$  are enough to determine the knot type, we talk of a simple minimal knot. Such a knot is given by three integers $N<p,q$; denoted by $K(N,p,q)$, it can be parametrized in the cylinder as $e^{i\theta}\mapsto (e^{Ni\theta},\sin q\theta,\cos p\theta)$. From this expression stems a natural representation of $K(N,p,q)$ as an $N$-braid. In this paper, we give a formula for its writhe number, i.e. the signed number of crossing points of this braid and derive topological consequences. We also show that if $q$ and $p$ are not mutually prime, $K(N,p,q)$ is periodic.
Simple minimal knots are a generalization of torus knots.
 \end{abstract}

\section{Introduction}
\subsection{Simple minimal knots}
In the present paper we pursue the investigation of simple minimal knots that we began in [So-Vi]. We recall that a simple minimal knot $K(N,p,q)$ is the intersection of the $3$-sphere $\mathbb{S}^3$ with $\phi(\mathbb{D})$ where 
$\mathbb{D}$ is the unit disk in $\mathbb{C}$ and $\phi$ is a map of the form 
$$\phi_{N,p,q}:\mathbb{D}\longrightarrow \mathbb{C}^2$$
$$\phi_{N,p,q}:z\mapsto (z^N, a(z^q-\bar{z}^q)+b(z^p+\bar{z}^p))$$
for some complex numbers $a$, $b$ and some integers $N<q,p$, with $(N,q)=(N,p)=1$.\\ 
We can write such a knot $K$ in $\mathbb{R}^3$ by the following expressions
$$x=2\cos Nt -\cos Nt\sin qt$$
$$y=2\sin Nt -\sin Nt\sin qt$$
$$z=\cos p t$$	
Note that if $p=q$, we get the $(N,q)$ torus knot.\\ Our motivation for studying these knots came from Riemannian geometry. 
We consider a minimal, i.e. conformal harmonic, map $F:\mathbb{D}\longrightarrow \mathbb{R}^4$ where 
$\mathbb{D}$ is the disk in $\mathbb{C}$; we assume that $dF(0)=0$, that is $F$ has a branch point at $0$.
If moreover that $F$ is a topological embedding, we can replicate Milnor's construction of algebraic knots and take the intersection of $F(\mathbb{D})$ with a small sphere centered at $F(0)$: we get a knot which we call a {\it minimal knot}. Simple minimal knots are in some sense the building blocks of minimal knots.
\section{The periodic case}
\begin{thm}
Let $d=gcd(p,q)$. If $d\neq 1$, then $K$ is a periodic knot.
\end{thm}
PROOF. We remind the reader ([So-Vi]) of the existence of a function $r:\mathbb{S}^1\longrightarrow \mathbb{R}_+^*$ such that 
$K=\{\phi_{N,p,q}(r(\theta)e^{i\theta}) / e^{i\theta}\in\mathbb{S}^1\}$. We notice that, for every $r>0$, 
$$\|\phi_{N,p,q} (r e^{i(\theta+\frac{2\pi}{d})})\|=\|\phi_{N,p,q} (r e^{i\theta})\|.$$
It follows that $r(\theta+\frac{2\pi}{d})=r(\theta)$. \\
We define an isometry $\Psi$ of $\mathbb{S}^3$ by setting for $(z_1,z_2)\in \mathbb{S}^3$,
$$\Psi(z_1,z_2)=(e^{\frac{2iN\pi}{d}}z_1,z_2).$$
The set $I$ of fixed points of $\Psi$ in $\mathbb{S}^3$ is $\{(0,e^{i\alpha}) / \alpha\in [0,2\pi]\}$ and verifies 
$lk(K,I)=\pm 1$; it is  diffeomorphic to $\mathbb{S}^1$ and disjoint from $K$. On the other hand, since $d$ divides $p$ and $q$ and $r(\theta+\frac{2\pi}{d})=r(\theta)$, we have 
$$\Psi(K)=K.$$
Thus $K$ is a periodic knot (for the definition of a periodic knot, see [Burde-Zieschang] p.266).
\section{The writhe number}
Simple minimal knots have a very natural braid representation; we will recall it and give a formula for its writhe number.
\subsection{The minimal braid}
We derive from the map $\phi_{N,p,q}$ above a braid $B$ in $Cyl=\mathbb{S}^1\times \mathbb{R}^2$ which naturally represents $K$: in [So-Vi] we called it the {\it minimal braid}. Its $N$ strands are given by the maps:
$$B_k:[\eta, 1+\eta]\longrightarrow \mathbb{R}^2$$
$$B_k:t\mapsto (\sin (\frac{2\pi}{N}q(t+k)),\cos (\frac{2\pi}{N}p(t+k)))$$
where $k=0,...,N-1$, $\eta$ is a very small positive irrational number (we introduce $\eta$ in order to avoid crossing points at the endpoints of the interval).\\
\\
REMARK. These knots can also be defined with a ${\it phase}$: the term $$\cos (\frac{2\pi}{N}p(t+k))$$ can be replaced by $$\cos (\frac{2\pi}{N}p(t+k+\epsilon)).$$ For technical reasons, we will use phases in our proofs below but we recall that the phase does not change the isotopy type of the braid up to mirror symmetry (see [So-Vi]). 
\subsection{The writhe number}
A crossing point $Q$ of the braid is the data of two different integers, $k$, $l$ with $0\leq k< l\leq N-1$ and a 
number $t\in [\eta, 1+\eta]$ such that
$$\sin (\frac{2\pi}{N}q(t+k))=\sin (\frac{2\pi}{N}q(t+l)).$$
There are $(N-1)q$ crossing points, as in the case of the torus knot where $q=p$; each crossing point $Q$ is assigned a sign $S(Q)\in\{-1,+1\}$ as below (in the case of the torus knot, they are all positive).\\
\begin{figure} 
 \picta  
\caption[fig2]{Sign of a crossing point }
\end{figure}
The writhe number $w(B)$ is the signed number of the crossing points, $$w(B)=\sum_{Q\mbox{\ crossing point}}S(Q).$$ It is related to the signed number of double points one gets when one desingularizes the branched disk $\phi(\mathbb{D})$ (see Introduction) into an immersion.

\subsection{Formulae for the writhe number}
\begin{prop}
Put $d=\gcd(p,q)$; let $\tilde{p}=\frac{p}{d}$ and $\tilde{q}=\frac{q}{d}$. \\
1) If $\tilde{p}$ and $\tilde{q}$ have different parities,
$$w(B)=0.$$
2) If $\tilde{p}$ and $\tilde{q}$ have the same parity, then up to sign,
$$w(B)=d\sum_{s=1}^{N-1}(-1)^{[\frac{ps}{N}]}(-1)^{[\frac{qs}{N}]}$$
where $[\ ]$ denotes the integral part.
\end{prop}

\begin{cor}
$$|w(B)|\leq d(N-1)\ \ \ \ \ \  (*)$$
The equality is attained in (*) if and only if $2N$ divides $p+q$ or $p-q$. Otherwise 
$$|w(B)|\leq d(N-3)\ \ \ \ \ \ \ \ $$
\end{cor}
By contrast we remind the reader that a $(N,q)$ torus knot has a writhe number of $(N-1)q$, since all the crossing points have a positive sign. \\
We can plug our estimates into the Franks-Williams-Morton inequalities for the HOMFLY polynomial ([F-W],[Mo]):
\begin{cor} Let $P_{min}$ (resp. $P_{max}$) be the minimal (resp. maximal) degree of the HOMFLY polynomial of $K(N,p,q)$. Then
$$-(d+1)(N-1)\leq P_{min}\ \ \mbox{and}\ \ \ P_{max}\leq (d+1)(N-1).$$
\end{cor}
If $p=q=d$ and $K$ is a torus knot,  both inequalities are equalities. On the other hand, if $p$ and $q$ are mutually prime, we get $$-2(N-1)\leq P_{min}\leq P_{max}\leq 2(N-1).$$
REMARK. A word of caution about Corollary 2: it cannot help us determine if a given knot is isotopic to a simple minimal knot. Indeed, $N$ is the number of strands of a specific braid representative of the knot $K$ but there is no reason why $K$ could not be represented by braids of smaller index. For instance, we saw in [So-Vi] that the knot $K(5,22,6)$ is isotopic to the knot $7_7$ (notation of the Rolfsen table ([Ro])) which has braid index $4$.
\section{Proofs}
We start with Prop. 1. 
1) We proved in [So-Vi] that if $p$ and $q$ have different parities, then $w(B)=0$; a similar proof works here. We use the notations of the introduction and compute
$$\Phi_{N,p,q}(re^{i(\theta+\frac{\pi}{d})})=(z^Ne^{iN\frac{\pi}{d}},
(-1)^{\tilde{q}}a(z^q-\bar{z}^q)+b(-1)^{\tilde{p}}(z^p+\bar{z}^p)).$$
Thus the transformation $H$ of $\mathbb{R}^4$ given below and the transformation that $H$ induces on $\mathbb{S}^3$ preserves $w(B)$:
$$H:\mathbb{C}\times\mathbb{R}\times\mathbb{R}\longrightarrow\mathbb{C}\times\mathbb{R}\times\mathbb{R}.$$
$$H:(w,x,y)\mapsto(we^{iN\frac{\pi}{d}},(-1)^{\tilde{q}}x,(-1)^{\tilde{p}}y)$$ On the other hand, if $\tilde{q}$ and $\tilde{p}$ have different parities, $H$ reverses the orientation. It follows that $w(B)=0$.\\
\\
We now assume that $\tilde{p}$ and $\tilde{q}$ have the same parity; since they are mutually prime, this means that 
$\tilde{p}$ and $\tilde{q}$ are both odd. We also derive that $p$ and $q$ have the same parity as well.\\
We give a proof of Prop. 1 2) based on Lemmas 1 to 6 below. \\
We have established in [So-Vi] that  $t\in[\eta, 1+\eta]$ is a crossing point between the $k$-th and $l$-th strands of $B$ if and only if there exists an integer $m$ such that $$t=-\frac{k+l}{2} + \frac{N}{4q}(2m+1)\ \ \ \ \ (0)$$
This prompts us to investigate the set $A$ of integers $m$ satisfying 
$$\eta<-\frac{k+l}{2} + \frac{N}{4q}(2m+1)\leq 1+\eta\ \ \ \ (1)$$ for some strands $k,l$. The largest (resp. smallest) value for $k+l$ is $1$, for $\{k,l\}=\{0,1\}$ (resp. $2N-3$ for $\{k,l\}=\{N-1,N-2\}$); hence $m$ belongs to $A$ if and only if
$$\eta+\frac{1}{2}\leq  \frac{N}{4q}(2m+1)\leq N-\frac{1}{2}+\eta\ \ \ \ (2)$$
\begin{lem}
Let $m$ be an integer belonging to $A$.\\
a) Then $m+q$ and $m-q$ cannot both belong to $A$\\
b) Assume that neither $m+q$ nor $m-q$ belongs to $A$ and let $k,l$ be the integers appearing with $m$ in (2). Then $k+l$ is one of the following three integers: $N-2$, $N-1$ or $N$. 
\end{lem}
PROOF. a) We see from 1) that the expression $\frac{N}{4q}(2m+1)$ belongs to an interval of amplitude $N-1$. On the other hand, we compute  
$$\frac{N}{4q}(2(m+q)+1)-\frac{N}{4q}(2(m-q)+1)=N$$
hence $\frac{N}{4q}(2(m+q)+1)$ and $\frac{N}{4q}(2(m-q)+1)=N$ cannot both satisfy a).\\
b) Assume that $m\in A$ but $m-q\notin A$ and $m+q\notin A$. We derive from (2)
$$\frac{N}{4q}(2(m-q)+1)<\eta+\frac{1}{2}\mbox{\ and\ }\frac{N}{4q}(2(m+q)+1)>N-\frac{1}{2}+\eta$$
hence 
$$\eta-\frac{1}{2}+\frac{N}{2}<\frac{N}{4q}(2m+1)<\eta+\frac{1}{2}+\frac{N}{2}.$$
Using the inequality (2), this yields
$$\frac{N}{2}-\frac{3}{2}<\frac{k+l}{2}<\frac{N}{2}+\frac{1}{2}$$
which translates into $N-2\leq k+l\leq N$, which proves b).\\
\\
Next, for a given $m$, we look for all possible values of $k,l$. Following Lemma 1 we  let $$A_0=\{m\in A\mbox{\ such that\ } m+q\notin A\mbox{\ and\ }m-q\notin A\}$$ 
$$A_1=\{m\in A\mbox{\ such that\ } m+q\in A\}.$$ 
Note that $A_0\cup A_1$ is included in $A$ but not equal to $A$; more precisely $A$ is the disjoint union $A=A_0\cup A_1\cup\{m+q\slash m\in A_1\}$. We have
\begin{lem}
Let $m$ be an integer. \\
i) If there exists an integer $S$ such that 
$$\eta\leq -\frac{S}{2}+\frac{N}{4q}(2m+1)\leq 1+\eta \ \ \ \ (3)$$
then there exist exactly two integers which verify (3). If we denote by $S_1(m)$ the smallest one of the two,  then the other one, $S_2(m)$, verifies $$S_2(m)=S_1(m)+1.$$
ii) If $m\in A_0$, then $S_1(m)=N-2$ or $S_1(m)=N-1$\\
iii) If $m\in A_1$, then $S_1(m+q)=S_1(m)-N$.
\end{lem}
REMARK. It follows that if $N=3$, $A_1$ is empty.\\
PROOF. i) If such integers exist for a given $m$, let $S_1(m)$ be the smallest one verifying (3). Then if $S$ is any integer such that $S\geq S_1(m)+2$, we derive from (3) that $$-\frac{S}{2}+\frac{N}{4q}(2m+1)\leq -1 -\frac{S_1(m)}{2}+\frac{N}{4q}(2m+1)< \eta$$
(the last inequality is strict because $\eta$ is irrational), hence $S$ does not satisfy (3). On the other hand, if $S_1(m)+1$ did not verify (3), a quick computation would show that $S_1(m)-1$ would  verify (3), which would contradict $S_1(m)$ being the smallest integer with that property. This proves i).\\
ii) If $m\in A_0$ then it verifies the assumptions of Lemma 1 b) and the possible values for $S$ are $N-2$, $N-1$ and $N$. It follows from Lemma 2 i) that $S_1(m)$ is $N-2$ or $N-1$.\\
To prove iii) it is enough to notice that 
$$-\frac{S}{2}+\frac{N}{4q}(2(m+q)+1)=-\frac{S-N}{2}+\frac{N}{4q}(2m+1).$$
Lemma 2 is now proved.\\
\\
We now introduce the phase $\epsilon$ which we mentioned but disregarded in \S 3.1: namely, we replace the $B_k$'s by  
$$B_k:t\mapsto (\sin (\frac{2\pi}{N}q(t+k)),\cos (\frac{2\pi}{N}p(t+k))+\epsilon).$$
To discuss signs of crossing points, we use the\\
NOTATION. If $r$ is a real number, $$\sigma(r)=(-1)^{[r]}\ \ \ \ \ \ \ \ \ (4)$$
We recall  the formula from [So-Vi] for the sign $S(m,k,l)$ of a crossing point given by (0)
$$S(m,k,l)=(-1)^{m}\sigma(p\frac{m}{q}+\epsilon)\sigma(p\frac{k-l}{N})\sigma(q\frac{k-l}{N})$$
To compute $w(B)$, we need to add the $S(m,k,l)$'s for all the $m,k,l$; we will group them according
$$w(B)=\sum_{m\in A_0} S(m,k,l)+ \sum_{m\in A_1} S(m,k,l)+\sum_{\{m+q\slash m\in A_1\}} S(m,k,l)\ \ \ \ \ \ (5)$$
We write the first term in (4) as $\sum_{m\in A_0}(-1)^{m}\sigma(p\frac{m}{q}+\epsilon)s_0(m)$ where 
$$s_0(m)=\sum_{k+l=S_1(m)}^{k+l=S_1(m)+1} \sigma(p\frac{k-l}{N})\sigma(q\frac{k-l}{N})\ \ \ \ \ \ \ \ \ \ \ (6)$$
We now want to write the sum of the second and third term in (5) as a sum over $A_1$ only; to this effect we
use the fact that $p$ and $q$ have the same parity (see above) to derive 
$$(-1)^{m+q}\sigma(p\frac{m+q}{q}+\epsilon)=(-1)^{p+q}(-1)^{m}\sigma(p\frac{m}{q}+\epsilon)
=(-1)^{m}\sigma(p\frac{m}{q}+\epsilon)\ \ \ (7).$$
Using (7), we can write the sum of the last two terms in (5) as 

$\sum_{m\in A_1}(-1)^{m}\sigma(p\frac{m}{q}+\epsilon)s_1(m)$ where 
$$s_1(m)=\sum_{k+l=S_1(m)}^{k+l=S_1(m)+1}\sigma(p\frac{k-l}{N})\sigma(q\frac{k-l}{N})
+\sum_{k+l=S_1(m)-N}^{k+l=S_1(m)-N+1} \sigma(p\frac{k-l}{N})\sigma(q\frac{k-l}{N})\ \ \ \ \ \ \ \ \ \ \ (8) $$
We put all this together to derive
$$w(B)=\sum_{m\in A_i\slash i=1,2}(-1)^{m}\sigma(p\frac{m}{q}+\epsilon )s_i(m)\ \ \ \ \ \ \ \ (9)$$
where $\sigma$, $s_0$ and $s_1$ are respectively defined in (4), (6) and (8).
\begin{lem} 
If $m$ belongs to $A_i$, $i=0,1$, we have $$s_i(m)=\sum_{s=1}^{N-1}\sigma(p\frac{s}{N})\sigma(q\frac{s}{N}) $$
\end{lem}
To prove Lemma 3, we will use the following easy identities ($n$ is an integer)
$$\sigma(p\frac{n}{N})\sigma(q\frac{n}{N})=\sigma(-p\frac{n}{N})\sigma(-q\frac{n}{N})=\sigma(p\frac{2N-n}{N})q(\sigma\frac{2N-n}{N})\ \ \ (10)$$
$$\sigma(p\frac{n}{N})\sigma(q\frac{n}{N})=(-1)^{p+q}\sigma(p\frac{N+n}{N})\sigma(q\frac{N+n}{N})
=\sigma(p\frac{N-n}{N})\sigma (q\frac{N-n}{N})\ \ \ \ (11)$$
Note that (10) is true for every set of integers $q,p,N$ whereas (11) is true because $p$ and $q$ have the same parity.\\
\\
We split the proof of Lemma 3 into three cases which are different but very similar: therefore we will omit repeating 
details. In all that follows we take
$k<l$, in other words $$k=min(k,l)$$

{\bf 1st case. $m\in A_0$ and $S_1(m)=N-2$}
$$s_0(m)=\sum_{k+l=N-2}^{k+l=N-1}\sigma(p\frac{k-l}{N})\sigma(q\frac{k-l}{N})$$
$$=\sum_{0\leq 2k<N-2}\sigma(p\frac{N-2-2k}{N})\sigma(q\frac{N-2-2k}{N})\ \ \ \ \ \ \ \ \ \ \ \ (12)$$
$$+\sum_{0\leq 2k<N-1}\sigma(p\frac{N-1-2k}{N})\sigma(q\frac{N-1-2k}{N})\ \ \ \ \ \ \ \ \ \ \ \ (13)$$
We rewrite (13) in terms of $2k-1$ instead of $2k$: we replace $N-1-2k$ by $N-2-(2k-1)$ and $0\leq 2k\leq2k-1$. We get 
$$(13)=
\sum_{-1\leq 2k-1<N-2}\sigma(p\frac{N-2-(2k-1)}{N})\sigma(q\frac{N-2-(2k-1)}{N})$$
Using this, we rewrite (12)+(13) as a single sum where we replace the variables $2k$ and $2k-1$ by a single 
variable $u$ such that $-1\leq u<N-2$. We get 
$$s_0(m)=\sum_{u=-1}^{N-3}\sigma(p\frac{N-2-u}{N})\sigma(q\frac{N-2-u}{N}).$$
We replace $u$ by $s=u+2$ and get 
$$s_0(m)=\sum_{v=1}^{N-1}\sigma(p\frac{N-s}{N})\sigma(q\frac{N-s}{N})=
\sum_{s=1}^{N-1}\sigma(p\frac{s}{N})\sigma(q\frac{s}{N})$$
the last identity following from (11).\\
\\
{\bf 2nd case. $m\in A_0$ and $S_1(m)=N-1$}\\
\\
The proof is identical to the 1st case except that, when we sum on $k+l=N$, we need to take  $2\leq 2k$ 
instead of $0\leq 2k$.
$$s_0(m)=\sum_{k+l=N-1}\sigma(p\frac{k-l}{N})\sigma(q\frac{k-l}{N})+
\sum_{k+l=N}\sigma(p\frac{k-l}{N})\sigma(q\frac{k-l}{N})$$
$$=\sum_{0\leq 2k<N-1}\sigma(p\frac{N-1-2k}{N})\sigma(q\frac{N-1-2k}{N})$$
$$+\sum_{1\leq 2k-1<N-1}\sigma(p\frac{N-1-(2k-1)}{N})\sigma (q\frac{N-1-(2k-1)}{N})$$ 
$$=\sum_{u=0}^{N-2}\sigma(p\frac{N-1-u}{N})\sigma(q\frac{N-1-u}{N})=
\sum_{s=1}^{N-1}\sigma (p\frac{N-s}{N})\sigma(q\frac{N-s}{N})$$
$$=\sum_{s=1}^{N-1}\sigma(p\frac{s}{N})\sigma(q\frac{s}{N})$$ 
{\bf 3rd case. $m\in A_1$}$\ \ \ \ \ \ \ \ \ \ \ \ \ \ \ s_1(m)=$
$$ \sum_{k+l=S_1(m)-N}\sigma(p\frac{k-l}{N})\sigma(q\frac{k-l}{N})
+\sum_{k+l=S_1(m)-N+1}\sigma(p\frac{k-l}{N})\sigma(q\frac{k-l}{N})\ (14)$$
$$+\sum_{k+l=S_1(m)}\sigma(p\frac{k-l}{N})\sigma(q\frac{k-l}{N})
+\sum_{k+l=S_1(m)+1}\sigma(p\frac{k-l}{N})\sigma(q\frac{k-l}{N})\ \ \ \ \ (15)$$
We start by dealing with (14) using the same method as in the previous two cases and write 
$$(14)=\sum_{0\leq 2k<S_1(m)-N}\sigma(p\frac{S_1(m)-N-2k}{N})\sigma(q\frac{S_1(m)-N-2k}{N})$$
$$+\sum_{-1\leq 2k-1<S_1(m)-N}\sigma(p\frac{S_1(m)-N-(2k-1)}{N})\sigma(q\frac{S_1(m)-N-(2k-1)}{N})$$

$$=\sum_{u=-1}^{S_1(m)-N-1}\sigma(p\frac{S_1(m)-N-u}{N})\sigma(q\frac{S_1(m)-N-u}{N})$$
$$=\sum_{s=N+1}^{S_1(m)+1}\sigma(p\frac{v}{N})\sigma(q\frac{s}{N})=\sum_{s=1}^{S_1(m)-N+1}\sigma(p\frac{s}{N})\sigma(q\frac{s}{N})$$

To treat (15), we proceed as in all the previous two cases, except that we index our sum with $l$ instead of $k$. If $k+l=S_1(m)$ (resp. $k+l=S_1(m)+1$), we have $2l>S_1(m)$ (resp. $2l-1>S_1(m)$). We derive
$$(15)=\sum_{S_1(m)<2l\leq 2N-2}\sigma(p\frac{S_1(m)-2l}{N})\sigma(q\frac{S_1(m)-2l}{N})$$
$$
+\sum_{S_1(m)<2l-1\leq 2N-3}\sigma(p\frac{S_1(m)-(2l-1)}{N})\sigma(q\frac{S_1(m)-(2l-1)}{N})$$
$$=\sum_{u=S_1(m)+1}^{2N-2}\sigma(p\frac{S_1(m)-u}{N})\sigma(q\frac{S_1(m)-u}{N})
=\sum_{s=S_1(m)-2(N-1)}^{-1}\sigma(p\frac{s}{N})\sigma(q\frac{s}{N})$$
$$=\sum_{s=2+S_1(m)-N}^{N-1}\sigma(p\frac{s}{N})\sigma(q\frac{s}{N}).$$
Adding (14) and (15) yield the third case of Lemma 3.\ \ \ \ \ \ $\square$ 

Putting Lemma 3 and formula (9) together, we derive
$$w(B)=\sum_{m\in A_i\slash i=1,2}(-1)^{m}\sigma(p\frac{m}{q}+\epsilon)\sum_{s=1}^{N-1}\sigma(p\frac{s}{N})\sigma(q\frac{s}{N}) \ \ \ \ \ \ \ \ \ (16)$$
We now focus on $\sum_{m\in A_i\slash i=1,2}(-1)^{m}\sigma(p\frac{m}{q}+\epsilon)$.
\begin{lem}
There exists an integer $m_0$ such that $A_0\cup A_1=[m_0, m_0+q-1]\cap\mathbb{N}$.
\end{lem}
PROOF. We go back to inequality (2) above and derive that $m$ belongs to $A$ if and only if
$$\frac{2q}{N}(\eta+\frac{1}{2})-\frac{1}{2}\leq m\leq 2q+\frac{2q}{N}(\eta-\frac{1}{2})-\frac{1}{2};$$
this defines an interval of amplitude $2q-\frac{2q}{N}$. We derive the existence of two integers $m_0$, $C$ with $0<C<q$
such that $A=[m_0, m_0+q+C]\cap\mathbb{N}$. \\
With these notations, we see that, if $m_0\leq m$,\\
1) $m\in A_1$ if and only if $m+q\leq m_0+q+C$, hence $A_1=[m_0, m_0+C]\cap\mathbb{N}$\\
2) $m\in A_0$ if and only if $m-q<m_0$ and $m+q>m_0+q+C$. In other words, $A_0=[m_0+C+1, m_0+q-1]\cap\mathbb{N}$.\\
The lemma follows from 1) and 2) put together.\\
\\
We derive from Lemma 4 that $$\sum_{m\in A_i\slash i=1,2}(-1)^{m}\sigma(p\frac{m}{q}+\epsilon)=
\sum_{u=1}^q\sigma(\frac{p(m_0-1+u)}{q}+\epsilon)(-1)^{m_0-1+u}$$
$$=
\sum_{u=1}^q\sigma(\frac{pu}{q}+\frac{p(m_0-1)}{q}+\epsilon+m_0-1)(-1)^{u}=\sum_{u=1}^q\sigma(\frac{pu}{q}+\psi)(-1)^u$$
for some irrational number $\psi$. Thus  we can rewrite (16) as 
$$w(B)=(\sum_{u=1}^q\sigma(\frac{pu}{q}+\psi)(-1)^u)(\sum_{s=1}^{N-1}\sigma(p\frac{s}{N})\sigma(q\frac{s}{N}))\ \ \ \ \ \  (17)$$
\begin{lem}  $\sum_{u=1}^q\sigma(\frac{pu}{q}+\psi)(-1)^u=
d\sum_{u=1}^{\tilde{q}}\sigma(\frac{\tilde{p}u}{\tilde{q}}+\psi)(-1)^u.$
\end{lem}
PROOF. $\sum_{u=1}^q\sigma(\frac{pu}{q}+\psi)(-1)^u=
\sum_{u=1}^{q}\sigma(\frac{\tilde{p}u}{\tilde{q}}+\psi)(-1)^u$
$$=\sum_{a=0}^{d-1}\sum_{v=1}^{\tilde{q}}\sigma(\frac{\tilde{p}(v+a\tilde{q})}{\tilde{q}}+\psi)(-1)^{v+a\tilde{q}}
=(\sum_{a=0}^{d-1}(-1)^{a(\tilde{p}+\tilde{q})})
\sum_{u=1}^{\tilde{q}}\sigma(\frac{\tilde{p}u}{\tilde{q}}+\psi)(-1)^u$$
$$=d\sum_{u=1}^{\tilde{q}}\sigma(\frac{\tilde{p}u}{\tilde{q}}+\psi)(-1)^u$$
since $\tilde{p}+\tilde{q}$ is even.

\begin{lem} Up to sign $\sum_{u=1}^{\tilde{q}}\sigma(\frac{\tilde{p}u}{\tilde{q}}+\psi)(-1)^u=1$
\end{lem}
PROOF. For every $u$, we have $\sigma(\frac{\tilde{p}(u+\tilde{q})}{\tilde{q}}+\psi)(-1)^{u+\tilde{q}}=\sigma(\frac{\tilde{p}u}{\tilde{q}}+\psi)(-1)^u$, hence
$$\sum_{u=1}^{\tilde{q}}\sigma(\frac{\tilde{p}u}{\tilde{q}}+\psi)(-1)^u=
\frac{1}{2}\sum_{u=1}^{2\tilde{q}}\sigma(\frac{\tilde{p}u}{\tilde{q}}+\psi)(-1)^u=
\frac{1}{2}\sum_{u=1}^{2\tilde{q}}\sigma(\frac{\tilde{p}u}{\tilde{q}}+\psi)(-1)^{\tilde{p}u}
;$$
this last identity comes from the fact that $\tilde{p}$ is odd.\\
Since $\tilde{p}$ and $2\tilde{q}$ are mutually prime, the map
$$\mathbb{Z}/2\tilde{q}\mathbb{Z}\longrightarrow \mathbb{Z}/2\tilde{q}\mathbb{Z}$$
$$u\mapsto \tilde{p}u$$
is a bijection which preserves the parity - by parity of an element of $\mathbb{Z}/2\tilde{q}\mathbb{Z}$, we mean the parity of its representative in $[1, 2\tilde{q}]$. Hence we have 
$$\frac{1}{2}\sum_{u=1}^{2\tilde{q}}\sigma(\frac{\tilde{p}u}{\tilde{q}}+\psi)(-1)^{\tilde{p}u}
=\frac{1}{2}\sum_{u=1}^{2\tilde{q}}\sigma\frac{u}{\tilde{q}}+\psi)(-1)^u=
\sum_{u=1}^{\tilde{q}}\sigma(\frac{u}{\tilde{q}}+\psi)(-1)^u
$$
We let $r=\psi-[\psi]$.
 Since the statement of the Lemma is up to sign, we can assume that $[\psi]$ is even, so for every $u$,
$\sigma(\frac{u}{\tilde{q}}+\psi)=\sigma(\frac{u}{\tilde{q}}+r)$. For every $u$ with $1\leq u\leq\tilde{q}$, we have $0<\frac{u}{\tilde{q}}+r<2$ since $r$ is irrational. It follows that 
$$\sigma(\frac{u}{\tilde{q}}+r)=1 
\Leftrightarrow
 \frac{u}{\tilde{q}}+r<1\Leftrightarrow u<\tilde{q}(1-r)\Leftrightarrow 1\leq u\leq [\tilde{q}(1-r)]$$
Similarly
$$\sigma(\frac{u}{\tilde{q}}+r)=-1\Leftrightarrow [\tilde{q}(1-r)]+1\leq u\leq\tilde{q}$$
The reader can now finish the proof of Lemma 6 him/herself, treating separately the cases when 
$[\tilde{q}(1-r)]$ is even and when $[\tilde{q}(1-r)]$ is odd. \\ 
\\
This concludes the proof of Prop. 1.\\
\\
PROOF OF COROLLARY 1. When we look at the formula for $w(B)$ in Prop. 1, we realize that equality is attained in (*) if and only if for every $s\in\{1,...,N-1\}$, $\sigma(p\frac{s}{N})\sigma(q\frac{s}{N})$ have the same sign.\\
\\
Suppose that $2N|p+q$ or $2N|p-q$; then there exists $\tau\in\{-1,+1\}$ and $\mu\in\mathbb{Z}$ such that 
$$p+\tau q= 2N\mu$$
and $p=2N\mu-\tau q$. It follows that, for every integer $s$,  
$$\sigma(\frac{ps}{N})\sigma(\frac{qs}{N})=(-1)^{2N\mu s+\tau}\sigma(\frac{ps}{N})\sigma(\frac{qs}{N})=(-1)^{\tau}$$
(we remind the reader that if $x$ is not an integer, $[-x]=-[x]-1$, hence $\sigma(-x)=-\sigma(x)$) hence the $\sigma(p\frac{s}{N})\sigma(q\frac{s}{N})$'s have all the same sign.\\
\\
The converse will follow from
\begin{lem}
i) If $\forall s\in\{1,...,N-1\}$, $\sigma(\frac{qs}{N})=\sigma(\frac{ps}{N})$, then $2N$ divides $p-q$\\
ii) If $\forall s\in\{1,...,N-1\}$, $\sigma(\frac{qs}{N})=-\sigma(\frac{ps}{N})$, then $2N$ divides $p+q$.
\end{lem}
PROOF. We assume i) of Lemma 7 and we do the Euclidean division of $p$ and $q$ by $N$: we get $$p=Na_1+R_1,\   q=Na_2+R_2\ \ \mbox{with}\ \ 1\leq R_1, R_2\leq N-1\ \ \ \ \ (18)$$
(the left inequality in (18) comes from the fact that $N$ does not divide $p$ or $q$).\\
If we write assumption i) for $s=1$, we get 
$$\sigma(a_1+\frac{R_1}{N})=\sigma(a_2+\frac{R_2}{N})$$
which implies that $$(-1)^{a_1}=(-1)^{a_2}\ \ \ \ \ \ \ \ \ (19)$$
that is, $a_1$ and $a_2$ have the same parity.\\  
\begin{lem}
Under the assumptions of Lemma 7 i), for every $s=1,...,N-1$
$$[\frac{R_1s}{N}]=[\frac{R_2s}{N}]\ \ \ \ (20)$$
\end{lem}
PROOF. Since $(N,p)=(N,q)=1$, it follows that for $i=1,2$, $(N,R_i)=1$. We keep this in mind and prove Lemma 8 by induction on $s$.\\
If $s=1$, (20) follows immediately from (18). We now assume that (20) is true for $s$: there exists an integer $M$ such that for $i=1,2$,
$$MN\leq sR_i<(M+1)N.$$
Since $(N,R_i)=1$, this implies
$$MN+1\leq sR_i\leq MN+N-1\ \ \ \ i=1,2$$
hence
$$M+\frac{R_i+1}{N}\leq \frac{R_i(s+1)}{N}\leq M+1+\frac{R_i-1}{N}.$$
Since $\frac{R_i+1}{N}<2$ and $\frac{R_i-1}{N}<1$, we derive that $$[\frac{R_i(s+1)}{N}]\in\{M,M+1\}\ \ \ \ \ \ \ \ \ \ \ (21)$$
On the other hand, we have $(s+1)p=(s+1)a_1N+(s+1)R_1$ and $(s+1)q=(s+1)a_2N+(s+1)R_2$; since $(s+1)a_1N$ and $(s+1)a_2N$ have the same parity (19), i) implies that 
$$\sigma(\frac{(s+1)R_1}{N})=\sigma(\frac{(s+1)R_2}{N})\ \ \ \ \ \ \ \ \ \ (22).$$
Putting (20) et (22) together yields $[\frac{(s+1)R_1}{N}]=[\frac{(s+1)R_2}{N}]$.$\ \ \ \ \ \ \ \ \ \square$.\\
\\
Going back to the proof of Lemma 7i), we write Lemma 8 for $s=N-1$ and we put $$E=[\frac{(N-1)R_1}{N}]=[\frac{(N-1)R_2}{N}].$$ For $i=1,2$, we have
$$EN\leq (N-1)R_i \leq EN+N-1$$
We derive that $(N-1)|R_1-R_2|\leq N-1$ which implies that $R_1=R_2$ or 
$|R_1-R_2|=1$. But we see from (18) that
$$R_1-R_2=p-q+N(a_1-a_2)$$
which is an even number since $p-q$ is even by assumption and $a_1-a_2$ is even by (19).
 hence $R_1=R_2$ and 
$$p-q=N(a_1-a_2)\ \ \ \ \ \ \ \ \ \ \ \ \ (23)$$
Since $a_1-a_2$ is even number by (19), this concludes the proof of Lemma 7 i).\\
\\
Suppose now that $p,q,N$ verify ii) of Lemma 7 and let $T$ be a positive integer such that $2NT-q>0$. We get, for $s=1,...,N-1$
$$\sigma((2NT-q)\frac{s}{N})=-\sigma(\frac{qs}{N})=\sigma(\frac{ps}{N}).$$
We let $2NT-q$ play the role which $q$ played above and we apply Lemma 7 i) to the integers $2NT-q,p,N$; we derive
$$2N|(2NT-q)-p.$$
It follows that $2N$ divides $p+q$ and this concludes the proof of Lemma 7 ii).\\
\\
To conclude the proof of Cor. 1, note that if the $\sigma(\frac{qs}{N})$'s and $\sigma(\frac{ps}{N})$'s
are not all of the same sign, then $|\sum_{s=1}^{N-1}\sigma(\frac{qs}{N})\sigma(\frac{ps}{N})|\leq N-3$.


\begin{thebibliography}{so-vi}
\bibitem[F-W]{f-w} Franks J., Williams R., {\it Braids and the Jones-Conway polynomial}, Trans. Amer. Math. Soc. 303 (1987) 97-108\\
\bibitem[Mo]{mo} Morton H.R., {\it Seifert circles and knot polynomials}, Math. Proc. Cambridge Philos. Soc. 99 (1986) 107-109\\
\bibitem[Ro]{Ro} D. Rolfsen, {\it Knots and Links}, Publish or Perish, Houston, 1990.
\bibitem[So-Vi]{So-Vi} M. Soret, M. Ville, {\it Singularity knots of minimal surfaces in $\mathbb{R}^4$} Journal of Knot theory and its ramifications, 20(04) 2011.\\
\end{thebibliography}
\end{document}